\documentclass[10pt,twosided]{article}
\usepackage[tbtags]{amsmath}
\allowdisplaybreaks[4]
\usepackage{color}

\usepackage{amssymb,color}

\definecolor{c20}{rgb}{0.,0.7,0.}
\definecolor{c30}{rgb}{0.,0.,1.}
\definecolor{c40}{rgb}{1,0.1,0.7}
\definecolor{c50}{rgb}{1,0,0}
\definecolor{c60}{rgb}{1,0.9,0.1}

\def\cE#1{\textcolor{black}{#1}}
\def\cE#1{#1}
\def\cJ#1{\textcolor{c50}{#1}}
\def\cJ#1{#1}
\def\pE#1{\textcolor{green}{#1}}
\def\pE#1{#1}

\def\dE#1{#1}

\usepackage{amssymb,color}

\newcommand{\kb}[1]{\boldsymbol{#1}}
\newcommand{\vk}[1]{\kb{#1}}
\def\kal#1{{\cal{ #1}}}

\newcommand{\ve}{\varepsilon}

\newcommand{\abs}[1]{\lvert #1 \rvert}

\newcommand{\E}[1]{\mathbb{E}\left(#1\right)}
\newcommand{\VA}[1]{\mathbb{V}ar\left(#1\right)}
\newcommand{\CO}[1]{\mathbb{C}\mathrm{ov}\left(#1\right)}

\newcommand{\pk}[1]{\mathbb{P} \left \{ #1 \right\} }

\newcommand{\pb}[1]{\mbox{\rm$\vk{P}$}\Bigl \{#1 \Bigr \}}

\newcommand{\R}{\!I\!\!R}

\newcommand{\limit}[1]{\lim_{#1 \to   \infty}}

\newcommand{\BQN}{\begin{eqnarray}}
\newcommand{\EQN}{\end{eqnarray}}
\newcommand{\BQNY}{\begin{eqnarray*}}
\newcommand{\EQNY}{\end{eqnarray*}}

\newcommand{\BS}{\begin{sat}}
\newcommand{\ES}{\end{sat}}
\newcommand{\BT}{\begin{theo}}
\newcommand{\ET}{\end{theo}}
\newcommand{\BK}{\begin{korr}}
\newcommand{\EK}{\end{korr}}

\newcommand{\BD}{\begin{de}}
\newcommand{\ED}{\end{de}}

\newcommand{\BIT}{\begin{itemize}}
\newcommand{\EIT}{\end{itemize}}
\newcommand{\BDI}{\begin{description}}
\newcommand{\EDI}{\end{description}}

\newcommand{\BRM}{\begin{remarks}}
\newcommand{\ERM}{\end{remarks}}

\newcommand{\BEL}{\begin{lem}}
\newcommand{\EEL}{\end{lem}}

\newtheorem{theo}{Theorem}[section]
\newtheorem{sat}[theo]{Proposition}
\newtheorem{de}[theo]{Definition}
\newtheorem{lem}[theo]{Lemma}

\newtheorem{example}[theo]{Example}
\newtheorem{korr}[theo]{Corollary}
\newtheorem{remark}[theo]{Remark}
\newtheorem{remarks}[theo]{Remarks}
\newtheorem{prop}[theo]{Proposition}

\newcommand{\nelem}[1]{{Lemma \ref{#1}}}
\newcommand{\neprop}[1]{{Proposition \ref{#1}}}
\newcommand{\netheo}[1]{{Theorem \ref{#1}}}

\newcommand{\prooflem}[1]{\textsc{Proof of Lemma} \ref{#1}}

\newcommand{\COM}[1]{}
\newcommand{\US}[1]{\underset{#1}\sup}

\newcommand{\QED}{\hfill $\Box$}

\def\th{\Lambda}

\def\cA{\cE{\kal{C}}}

\def\XT{ \dE{\{X(t), t\in [0,T]\}}}
\def\XTI{ \dE{\{X(t), t\ge 0\}}}

\date{}

\def\rw{\rightarrow}
\def\IF{\infty}
\def\logsim{\overset{log}\sim}
\def\van{ \varepsilon}
\def\tO{t_0 }
\def\TU{T_u}
\def\eU{\ve_u}
\def\Co{\mathbb{C}\mathrm{ov}}

\def\PH{\mathcal{H}}
\def\LT{\left}
\def\RT{\right}
\newcommand{\OU}[1]{o\left(#1\right)}
\def\TT{\mathcal{T}}
\newcommand{\WW}[1]{\mathcal{W}\LT(#1\RT)}

\def\tO{t_0 }

\def\TU{\TT_u}
\def\eU{\ve_u}
\def\AA{{\bf A1-A4}}
\newcommand{\Lsup}[1]{\underset{#1}\limsup}
\def\const{\mathbb{\cJ{Q}}}

\def\pXT{p(X,\TT,\gamma,u)}

\newcommand{\tn}[1]{{\textcolor{black}{#1}}}

\def\tb#1{#1}

\begin{document}
\title{\bf \Large Tail Asymptotics of Supremum of Certain Gaussian Processes over Threshold Dependent Random Intervals}
\author{Krzysztof D\c{e}bicki\footnote{Mathematical Institute, University of Wroc\l aw, pl. Grunwaldzki 2/4, 50-384 Wroc\l aw, Poland},
Enkelejd Hashorva\footnote{
University of Lausanne, Faculty of Business and Economics (HEC Lausanne), 1015 Lausanne, Switzerland}, 
Lanpeng Ji$^\dagger$
}

 \maketitle
        \centerline{University of Wroc\l aw and University of Lausanne}
\vskip 0.41 cm
        \centerline{\today{}}

\vskip 1 cm

{\bf Abstract:} Let $\{X(t),t\ge0\}$ be a centered Gaussian
process and \pE{let} $\gamma$ be a non-negative constant. In this paper we
study the asymptotics of $\pk{\underset{t\in [0,\TT/u^\gamma]}\sup X(t)>u}$
as $u\rw\IF$, \cJ{with $\TT$} an
independent of $X$ non-negative random variable. As an
application, we derive the asymptotics of finite-time ruin
probability of time-changed fractional Brownian motion risk
processes.

\vskip 0.5 cm


\vskip 0.5 cm


\section{Introduction}

One of the seminal results in the extreme value theory of Gaussian processes concerns the asymptotic behavior
of
\begin{eqnarray}
\pk{\sup_{t\in[0,T]}X(t)>u}, \quad u\to\infty,
\label{eq.basic}
\end{eqnarray}
where $\{X(t),t\ge0\}$ is a centered Gaussian process with almost surely (a.s.) continuous sample paths,
variance function that attains its unique maximum at exactly one point $t_0\in[0,T]$
and $T>0$ is deterministic; see the classical monograph
Piterbarg (1996), Adler and Taylor (2007) for a complete survey on
distributional properties of extremes of Gaussian processes and related topics.

The aim of this paper is to analyze the asymptotics of
a counterpart of (\ref{eq.basic}) in the case of
\tn{
$T$ being an independent of $X$ non-negative random variable, possibly depending on
$u$, \cJ{namely}
 \BQN\label{eqXTT}
\pXT:=\pk{\underset{t\in [0,\TT/u^\gamma]}\sup X(t)>u}, \quad \text{ with } 
\gamma \ge 0
\EQN
as $u\to \IF$.
}

Let in the following $B_\alpha$ denote a fractional Brownian
motion (fBm) with Hurst index $\alpha/2\in(0,1]$, and let
$\{Y(t),t\ge0\}$ independent of $B_\alpha$ be a
non-negative,  non-decreasing stochastic process.\\
Apart from theoretical interest in properties of (\ref{eqXTT}),
the motivation to analyze \pE{it comes for numerous applications}.
 For instance the calculation of the finite-time ruin probability
of the time-changed fractional Brownian motion risk process $\psi_T$ defined by
\BQN\label{eqfbmY}
\psi_T(u)=\pk{\inf_{t\in[0,T]}(u+cY(t)-B_{\alpha}(Y(t)))<0}
\EQN
for some $c>0, T>0, u\ge0$ is possible using Monte Carlo simulations only when $Y$ is known.  In view of our findings,
for all large $u$ compact formulas for approximation of $\psi_T(u)$ can be given explicitly.\\
In the context of theoretical actuarial models, $u$ is the so-called {\it initial reserve}, $c$ models the premium income rate,
and $B_\alpha(t)$ represents the total amount up to time $t>0$ of aggregate claims (including fluctuations).
The justification for choosing fBm to model the aggregate claim process comes from Michna (1998),
where it is shown that the finite-time ruin probability $\psi_T(u)$, with $Y(t)=t$
is a good approximation of the finite-time ruin probability for the classical risk process with claims possessing {\it long range dependence} property.
 The role of the random process $Y$ is crucial in order to make such models adequate for applications. It is a substitute for the real time, where
 $Y(t)$ stands for the random business  time, which is consistent with the insurance practice where
both claims  and premiums may not be received immediately at time $t$ of the event, but at a later random time modeled by $Y(t)$.
Indeed, if $\{Y(t),t\ge0\}$ has additionally a.s. continuous sample paths,
then re-writing \eqref{eqfbmY} as
\BQN\label{eqfbmY:2}
\psi_T(u)
= \pk{\underset{t\in[0,\mathcal{T}/u ]}\sup \frac{B_{\alpha}(t)}{1+ct} >u^{1-\alpha/2}}
, \quad \text{ 
$\mathcal{T}= Y(T)$}
\EQN
we see that the asymptotic analysis of $\psi_T(u)$ reduces to the analysis of the \pE{asymptotic tail behaviour} of
supremum of a specific Gaussian process over a random interval.

Other branch of motivations to analyze (\ref{eqXTT}) stems from
recently studied problems in fluid queueing theory.
In particular, the tail asymptotics of the stationary buffer content
of a {\it hybrid fluid queue}, with input modelled by a superposition of integrated
alternating {\it on-off} process and a Gaussian process with stationary increments,
reduces (under some assumptions) to (\ref{eqXTT}) with suitably chosen random
$\TT$; \cJ{see} e.g., Zwart et al. (2005).

In the  case that $\{Y(t), t\ge 0\}$ is a Gamma process and $c=0$,
investigation of $\psi_T(u)$, as $u\to \IF$
reduces to the analysis of the tail asymptotics of the fractional Laplace motion; see Kozubowski et al. (2006),
Arendarczyk and D\c{e}bicki (2011).
For the related literature on time-changed models we refer to Fotopoulos and Luo (2011), who considered the case of Brownian motion ($\alpha=1$) 
and Wu and Wang (2012), where a model  based on the Cox risk process, which is a time-changed compound Poisson risk process, is considered.
For more applications 
in finance and insurance, see e.g., Geman et al. (2001) and Palmowski and Zwart (2007), among many others.

{\it Contribution}: a) In Theorem \ref{ThmT}, under some canonical asymptotic assumptions of the Gaussian process
$\{X(t),t\ge0\}$, see Section \ref{s.notation},
if $\TT$ has \pE{a} sufficiently heavy tail distribution
(which is manifested by {\it insensitivity} property of the tail distribution of $\TT$),
we derive exact asymptotics of $\pXT$. \\
b) Theorem \ref{ThmLog} deals with a more general class of
random variables $\TT$. Under a log-power tail asymptotic assumption (see \eqref{limL}), we obtain the logarithmic asymptotics of $\pXT$.
It appears that, depending on the interplay between heaviness of $\TT$
and local properties of the variance function of $X$ at $0$ we can distinguish
four scenarios, leading to four qualitatively different asymptotics.

Our novel result complement recent findings of Arendarczyk and D\c{e}bicki (2011, 2012) and Tan and Hashorva (2013),
where extremes over a random-time interval were analyzed for
stationary centered Gaussian processes and centered Gaussian processes with stationary increments respectively.
We refer to Zwart et al. (2005), D\c{e}bicki and van Uitert (2006), Palmowski and Zwart (2007)
for other related results.

 \pE{The} organization of the paper \pE{is as follows}: In Section 2 we introduce some notation and formulate the main assumptions imposed
to the Gaussian process $\XTI$.
Section 3 presents the main results.
In Section 4 we discuss the finite-time ruin probability of the time-changed fBm risk processes.
The proofs of all the results are relegated to Section 5.

\section{Preliminaries}\label{s.notation}

Throughout this paper  $\Psi(\cdot)$ denotes the survival function
of a standard Gaussian random variable $N(0,1)$. For two positive
functions $a_1(\cdot),a_2(\cdot)$ on $[0,\IF)$ we write
$a_1(u) \sim a_2(u)$ 
if $\lim_{u\to \IF} a_1(u)/a_2(u)=1$ and $a_1(u)\logsim a_2(u)$
if $\lim_{\cJ{u \to \IF}} \log(a_1(u))/\log(a_2(u))=1$. 
In the sequel  $\Gamma(\cdot)$ denotes  the Euler Gamma function.

Following, e.g., Foss et al. \pE{(2013)}
we say  that a non-negative random variable $\TT$ is $h-$insensitive if
$$\pk{\TT>u\pm h(u)}\sim \pk{\TT>u}$$
  as $u\to\infty$ for some function $h(\cdot)$. Our first main result in this paper is the exact asymptotic behaviour of
 $ \pXT$ defined in \eqref{eqXTT} for $\TT$ being $h-$insensitive. Two large  classes of \pE{distributions} for $\TT$ that satisfy the
 insensitivity criteria are:\\
  i) $\TT$ is {\it regularly varying} at $\infty$ with index $\lambda>0$, which means that
$\pk{\TT>u}=\mathcal{L}(u)u^{-\lambda}$, where $\mathcal{L}(\cdot)$ is slowly varying at $\infty$, i.e., for any $x>0$ we have
$\limit{u} \mathcal{L}(xu)/\mathcal{L}(u)=1$ \pE{(see e.g., Resnick (1987))};\\
ii)
\COM{  ifthat satisfy the
where $\TT$ is an independent of $X$ non-negative random variable with absolutely continuous \pE{distribution}.
In particular we analyze $\pXT$ for some classical classes of distributions of $\TT$.
We say that random variable $\TT$ is {\it regularly varying} at $\infty$ with index $\lambda\ge0$ if
$\pk{\TT>u}=\mathcal{L}(u)u^{-\lambda}$, where $\mathcal{L}(\cdot)$ is slowly varying at $\infty$, i.e., for any $x>0$ we have
$\limit{u} \mathcal{L}(xu)/\mathcal(u)=1$.
}
$\TT$ is {\it asymptotically Weibullian}, i.e.,
 \BQN\label{eqwei}
\pk{\TT> u}\sim  \mathcal{L}(u) u^{\delta}\exp(-Lu^p) 
\EQN
as $u\to\infty$, where $L, p>0$ and $\delta\in\R$ and $\mathcal{L}(\cdot)$ is  slowly varying  at $\IF$;
we abbreviate \eqref{eqwei} as $\TT\in\WW{p,L,\mathcal{L},\delta}$.\\

A significant relaxation of the Weibullian-tail assumption \eqref{eqwei}
for $\TT$ is  that it has asymptotically a {\it log-power} tail with coefficient $L>0$ and power $p>0$, i.e.,
\BQN\label{limL}
\underset{u\rightarrow\IF}\lim\frac{\log \pk{\TT>u}}{u^{p}}=-L.
\EQN

\COM{
The following result, due to
Piterbarg and Prisyazhnyuk (see e.g., Piterbarg (1996)), plays a crucial role in
further analysis.
\BT\label{ThmPP}
Let $\XT$ be a centered Gaussian process with standard deviation function $\sigma(\cdot)$ attaining its unique maximum point over $[0, T]$ at $\tO\in(0,T)$ with $\sigma(\tO)=1$. If further  {\bf A1-A3} are satisfied, then we have
\BQN
\pb{\sup_{t\in[0,T]}X(t)>u}&\sim & {\mathcal{C}}_{{{\alpha}},{{\beta}}} \th_{{{\alpha}}, {{\beta}}}(u) \Psi(u),
\quad u\to \infty,
\EQN
where
\BQNY
\cA_{{{\alpha}}, {{\beta}}}=\left\{
              \begin{array}{ll}
2\mathcal{H}_{{{\alpha}}}\Gamma(1/{{\beta}}+1)d^{1/{{\alpha}}}a^{-1/{{\beta}}}     , & \hbox{if } {{\alpha}}< {{\beta}} ,\\
\mathcal{P}_{{{\alpha}}}^{a/d}, & \hbox{if } {{\alpha}} ={{\beta}},\\
1, & \hbox{if } {{\alpha}} >{{\beta}},
              \end{array}
            \right.
\quad\mathrm{and}\quad
\th_{{{\alpha}},{{\beta}}}(u)=\left\{
              \begin{array}{ll}
u^{2/{{\alpha}}-2/{{\beta}}},& \hbox{if } {{\alpha}} <{{\beta}},\\
1, & \hbox{if } {{\alpha}} \ge{{\beta}}.
              \end{array}
            \right.
\EQNY
\ET

Theorem \ref{ThmPP} combined with the Borell-TIS inequality (e.g., Adler and Taylor (2007))
straightforwardly implies
Corollary \ref{K1}.
}

\def\PU{{\mathcal{C}}_{{{\alpha}},{{\beta}}} \th_{{{\alpha}}, {{\beta}}}(u) \Psi(u)}

As above, hereafter $B_\alpha$ stands for a fBm with Hurst index $\alpha/2\in(0,1]$, i.e.,
a centered Gaussian process with stationary increments, continuous sample paths a.s. and variance function
$\cJ{\sigma_{B_\alpha}^2(t)}=t^{\alpha}$.
Following Piterbarg (1996), we introduce two key constants, namely Pickands constant
$$\mathcal{H}_{\alpha}:=\lim_{S\rightarrow\infty} S^{-1} \E{ \exp\biggl(\sup_{t\in[0,S]}\Bigl(\sqrt{2}B_{\alpha}(t)-t^{\alpha}\Bigr)\biggr)}\cJ{\in(0,\IF)}$$
and Piterbarg 
constant
$$\mathcal{P}_{\alpha}^{R}:={\lim_{S \to \infty}}\E{ \exp\biggl(\sup_{t\in[-S,S]}\Bigl(  \sqrt{2} B_{\alpha}(t)-(1+R)t^{\alpha}\Bigr )\biggr)}\cJ{\in(0,\IF)}, \quad R\in (0,\IF).$$

In this paper we tacitly assume that $\XTI$ is a centered Gaussian process with a.s.
continuous and bounded sample paths and
$\sigma^2(t):=\VA{X(t)}$ that attains
its maximum over $[0,\IF)$ at the unique point $t=\tO>0$ with $\sigma(\tO)=1$.
Additionally we suppose that:\\
{\bf A1}.  There exist some positive constants $a,\beta$ such that
\BQN\label{eqA1}
\sigma(t)&=& 1-a\abs{t-\tO}^{\beta}+o(\abs{t-\tO}^{\beta}), \quad t\to \tO.
\EQN
\tn{
{\bf A2}.  There exist $d>0, \alpha\in {(0, 2]}$ such that
$${\CO{\frac{X(s)}{\sigma(s)}, \frac{X(t)}{\sigma(t)}}=1- d|t-s|^{\alpha}+o(|t-s|^{\alpha})}, \quad
s\rightarrow \tO, t\rightarrow \tO.$$
}

{\bf  A3}. There exist constants $\const>0$, $H >t_0$  and $r\in(0,2]$ such that, for all $s,t\in[0,H]$
 with $\abs{s-t}<1$
\BQN
\E{(X(t)-X(s))^{2}} &\leq & \const|t-s|^{r}.
\EQN


{\bf A4}. 
$\underset{t\rightarrow \infty}\limsup\ \sigma(t)<1$.

\bigskip
We conclude this section with  a preliminary result, which gives the exact asymptotics of the supremum of a Gaussian process over a deterministic time interval; see e.g., Piterbarg (1996).

\bigskip
\BT  \label{K1} If the assumptions {\bf A1-A4} are satisfied, then for any  $T>\tO$ we have
\BQN\label{eqXIF}
\pk{\US{t\in[0,\IF)}X(t)>u}
\sim \pk{\sup_{t\in[0,T]}X(t)>u} \sim \PU, \quad u\to \IF,
\EQN
where
\BQNY
\cA_{{{\alpha}}, {{\beta}}}=\left\{
              \begin{array}{ll}
2\mathcal{H}_{{{\alpha}}}\Gamma(1/{{\beta}}+1)d^{1/{{\alpha}}}a^{-1/{{\beta}}}     , & \hbox{if } {{\alpha}}< {{\beta}} ,\\
\mathcal{P}_{{{\alpha}}}^{a/d}, & \hbox{if } {{\alpha}} ={{\beta}},\\
1, & \hbox{if } {{\alpha}} >{{\beta}},
              \end{array}
            \right.
\quad\mathrm{and}\quad
\th_{{{\alpha}},{{\beta}}}(u)=\left\{
              \begin{array}{ll}
u^{2/{{\alpha}}-2/{{\beta}}},& \hbox{if } {{\alpha}} <{{\beta}},\\
1, & \hbox{if } {{\alpha}} \ge{{\beta}}.
              \end{array}
            \right.
\EQNY
\ET

{\bf Remark}: The exact asymptotics for the infinite-time interval case can be obtained by a direct application of the Borell-TIS inequality (e.g., Adler (1990) or Adler and Taylor (2007)) using further the result for the asymptotics of the supremum of $X$ over any finite-time interval $[0.T]$.

\section{Main Results}

In this section we present our main results.
\tb{We begin with the derivation of the}
exact asymptotic behaviour of $\pXT$ as $u\to\infty$, which is presented in Theorem \ref{ThmT}.
Then, under milder assumptions on $\TT$, we provide a complete study of the logarithmic asymptotics of
$\pXT$ as $u\to\infty$, see Theorem \ref{ThmLog}.

\BT\label{ThmT}
\cJ{Let $\{X(t), t\ge0\}$ be a centered Gaussian process  such that assumptions \AA \ are satisfied, and}
let $\TT$ be a non-negative random variable independent of the Gaussian process $X$.

i) If $\gamma=0$ and
$\pk{\TT\ge  \tO }>0$, then 
\BQN\label{MT}
\pXT
\sim  \PU  \pk{\mathcal{T} \ge  \tO }.
\EQN
ii) If $\gamma>0$
and $\mathcal{T}$ is $u^{1-2/(\gamma(1+\beta))}$-insensitive, then
\BQN\label{MT1}
\pXT
\sim  \PU  \pk{\mathcal{T}\ge  \tO u^\gamma }.
\EQN
\ET
The proof of Theorem \ref{ThmT} is given in Section \ref{s.proof.1}.

Theorem \ref{ThmT} complements recent results of
Arendarczyk and D\c{e}bicki (2011, 2012) and Tan and Hashorva (2013),  where the class of
centered stationary Gaussian processes and Gaussian processes with stationary increments
was analyzed.

As a straightforward consequence of Theorem \ref{ThmT} we obtain the following results.
\BK
\tn{Suppose that the assumptions of Theorem \ref{ThmT} are
satisfied and that $\gamma>0$.
}
\\
i) If $\TT$ is regularly varying at $\infty$ with index $\lambda>0$, then
\BQN\label{nase}
\pXT \sim  \PU  \pk{\mathcal{T}\ge  \tO u^\gamma }. 
\EQN
ii) If $\TT\in\mathcal{W}(p,L, \mathcal{L},\delta)$ with $p\in\left(0,\frac{2}{\gamma(1+\beta)}\right)$, then \eqref{nase} holds.
\COM{
\BQN
\pXT \sim  \PU  \pk{\mathcal{T}\ge  \tO u^\gamma }.
\EQN
}
\EK

Complementary to the above exact asymptotics, in the next theorem we derive
the logarithmic asymptotic behaviour of \eqref{eqXTT} for a class of log-power tailed random variables $\TT$.

Let us introduce the following notation
\BQN
\cE{\widehat \sigma(s):=\underset{t\in[0,s]}\sup \sigma(t), \quad
\widetilde \sigma_{L,\gamma}(s) :=\frac{1}{2\widehat \sigma^2(s)}+Ls^{2/\gamma}, \quad s\ge0}\label{not.1}
\EQN
and
\BQN
A_0=\inf\LT\{A: A=arg \underset{t\le t_0}\inf\LT( { \widetilde \sigma_{L,\gamma}(t)}\RT)\RT\}.\label{not.2}
\EQN

\BT\label{ThmLog} \cJ{Let $\{X(t), t\ge0\}$ be a centered Gaussian
process  such that assumptions \AA \ are satisfied} \tn{and let} $\TT$ be a
non-negative random variable independent of the Gaussian process
$X$
with asymptotically log-power tail with coefficient $L>0$ and power $p>0$.\\
i) If $\gamma p<2$, then
\BQN\label{eqLogl}
\lim_{u\to\infty}\frac{\log\pXT}{u^2} =-\frac{1}{2}.
\EQN

ii) If $\gamma p=2$, $A_0>0$, and on any compact subset of $(0,\infty)$ $\widehat \sigma(s)$ is differentiable and  
$\abs{\widehat \sigma'(s)}$ is uniformly bounded,
then
\BQN\label{eqLoge}
\lim_{u\to\infty}\frac{\log\pXT}{u^2}=
- \widetilde \sigma_{L,\gamma}(A_0).
\EQN

iii) If $\gamma p>2$ and $\sigma(0)>0,$ then
\BQN\label{eqLoggg}
\lim_{u\to\infty}\frac{\log\pXT}{u^2}=-\frac{1}{2\sigma^2(0)}.
\EQN

iv) If  $\gamma p>2$ and   $\sigma(t)=Dt^\eta(1+\OU{1})$ as $t\downarrow0$
for some positive constants $D$ and $\eta$, then
\BQN\label{eqLogD}
\lim_{u\to\infty}\frac{\log\pXT}{u^{2p(\eta\gamma +1)/(2\eta+p)}}
=-A_1,
\EQN
where $A_1=\frac{1}{2}
D^{-\frac{2p}{2\eta+p}}
(Lp/\eta)^{\frac{2\eta}{2\eta+p}}+L^{\frac{2\eta}{2\eta+p}}(\eta/(pD^2))^{\frac{p}{2\eta+p}}$.
\ET
A complete proof of Theorem \ref{ThmLog} is given in Section \ref{s.proof.2}.

\begin{remark}\label{rem1}
\tb{It follows straightforwardly from the proof of Theorem \ref{ThmLog} that
statements $i)$ and $iii)$ also hold if
$-\log \pk{\TT> u}= \mathcal{L}(u) u^p$ with some slowly varying function $\mathcal{L}(\cdot)$ at infinity.}

When $p \gamma=2$ we imposed the assumption that $A_0>0$. The special case $A_0=0$, which is also possible,
is much more involved, and will therefore \tb{be}
considered elsewhere.
\end{remark}

\COM{
Motivation (assumption):
\BQNY
\pk{\TT_u>u^{\delta}}<<\exp(-const u^2), \ \ \pk{\TT>u^{\delta+\gamma}}\overset{Logar}\sim \exp(-u^{(\delta+\gamma)p}),
\EQNY
which gives
$$
\delta> \frac{2-\gamma p}{p}.
$$
}

\section{Ruin probability of the time-changed fBm risk processes}

\tb{Consider an extension of the time-changed fBm risk process defined in the Introduction,
by allowing a power trend-function; i.e., let}
\BQN\label{ZZ}
Z(t):=u+c(Y(t))^\theta-B_{\alpha}(Y(t)),\ \ \ t\ge0,
\label{eqfbm1}
\EQN
where 
$ u\ge 0, c>0,\theta>\alpha/2$ and $\{Y(t), t\ge 0\}$ is a non-negative non-decreasing 
stochastic process being independent of $\{B_{\alpha}(t), t\ge 0\}$.
\tb{Clearly,
$\theta=1$
is a
choice leading to our risk model in the Introduction.
Related risk models were studied
for instance in D\c{e}bicki and Rolski (2002), where the finite-time ruin probability with the choice
$Y(t)\equiv t$ was analyzed, whereas the infinite-time ruin counterpart was considered in H\"{u}sler and  Piterbarg (1999), \pE{see also the recent
contribution Hashorva et al. (2013).}
}

As in the Introduction the finite-time ruin probability for the risk process
\eqref{ZZ} is given by
\BQNY
\psi_T(u)=\pk{\underset{t\in[0,T]}\inf \Bigl( u+c(Y(t))^\theta-B_{\alpha}(Y(t))\Bigr)<0}=
\pk{\underset{t\in[0,T]}\sup \Bigl( B_{\alpha}(Y(t))-c(Y(t))^\theta \Bigr)>u},
\EQNY
with $T>0$ and $u\ge 0$.
\subsection{Continuous time-process $Y$}
In this subsection, we apply the results of Theorems \ref{ThmT} and \ref{ThmLog} to the  analysis of the asymptotics of
the finite-time ruin probability of the time-changed fBm risk process \eqref{ZZ} as $u\rw\IF$, where the time process $Y$ has a.s. continuous sample paths. 

Before stating our results for this risk model we need to introduce the following notation
\begin{eqnarray*}
Q&:=&2^{\frac{1}{2}+\frac{1}{\alpha}}\sqrt{\pi}c^{\frac{2-\alpha}{2\theta}}\alpha^{\frac{\alpha-2-\theta}{2\theta}}
\theta^{\frac{2-\alpha}{\alpha}}(2\theta-\alpha)^{\frac{\theta\alpha-4\theta+2\alpha-\alpha^2}{2\theta\alpha}},\\
s_0&:=&\LT(\frac{\alpha}{c(2\theta-\alpha)}\RT)^{1/\theta}, \quad V_0:= \frac{2\theta-\alpha}{2\theta}s_0^{\alpha/2}= \frac{2\theta-\alpha}{2\theta}\LT(\frac{\alpha}{c(2\theta-\alpha)}\RT)^{\frac{\alpha}{2\theta}}.
\end{eqnarray*}

The main results are presented in Propositions \ref{Thmfbm1} and \ref{Thmfbm2}; 
their proofs are relegated to Section \ref{s.prop1_2}.

\begin{prop}\label{Thmfbm1}
Assume that $\theta>\alpha/2$, $c>0$ and $\{Y(t), t\ge 0\}$ has a.s. continuous sample paths.\\
i) If $Y(T)$ is regularly varying at $\IF$ with index $\lambda>0$, i.e.,
$\pk{Y(T)>u}=\mathcal{L}(u)u^{-\lambda}$, then
\BQN
\psi_T(u)\sim  Q\PH_\alpha s_0^{-\lambda}
u^{\frac{(2\theta-\alpha)(2-\alpha)-2\lambda\alpha}{2\theta\alpha}}\mathcal{L}(s_0 u^{\frac{1}{\theta}})\Psi
\LT(V_0^{-1}u^{\frac{2\theta-\alpha}{2\theta}}\RT). 
\EQN

ii) If $Y(T)\in\WW{p,L, \mathcal{L},\delta}$ with $p\in \left(0, \frac{2\theta-\alpha}{3}\right)$, $L>0$, and $\delta\in\R$,
then
\BQN
\psi_T(u)\sim  Q\PH_\alpha \mathcal{L}( s_0 u^{1/\theta})s_0^\delta
u^{\frac{(2\theta-\alpha)(2-\alpha)+2\delta\alpha}{2\theta\alpha}}\exp\LT(-L s_0^p u^{\frac{p}{\theta}}\RT)\Psi
\LT(V_0^{-1}u^{\frac{2\theta-\alpha}{2\theta}}\RT). 
\EQN
\end{prop}

\begin{prop}\label{Thmfbm2}
Under the setup of \neprop{Thmfbm1}, suppose further
\tb{that $Y(T)$ has asymptotically log-power tail with coefficient $L>0$ and power $p>0$.}

i) If $2\theta-\alpha>p$, then
\BQN
\lim_{u\to\infty}\frac{\log \psi_T(u)}{u^{\frac{2\theta-\alpha}{\theta}}}
= \frac{-1}{2V_0^{2}}.
\EQN

ii) If $2\theta-\alpha<p$, then
\BQN
\lim_{u\to\infty}\frac{\log \psi_T(u)}
                      { u^{\frac{2p}{\alpha+p}} }
=
-\LT(\frac{1}{2}\LT(\frac{2p}{\alpha}\RT)^{\frac{\alpha}{\alpha+p}}+\LT(\frac{\alpha}{2p}\RT)^{\frac{p}{\alpha+p}}\RT)L^{\frac{\alpha}{\alpha+p}}.
\EQN

iii) If $2\theta-\alpha=p$, then
\BQN
\lim_{u\to\infty}\frac{\log \psi_T(u)}
{u^{\frac{2\theta-\alpha}{\theta}}}
=
-\LT(\frac{(1+c A_0^\theta)^2}{2A_0^\alpha}+LA_0^{2\theta-\alpha}\RT), 
\EQN
where
$$
A_0=\left(\frac{c(\alpha-\theta)+\sqrt{c^2(\theta-\alpha)^2+2\alpha\LT(c^2(\theta-\alpha/2)+L(2\theta-\alpha)\RT)}}
{c^2(2\theta-\alpha)+2L(2\theta-\alpha)}\right)^{1/\theta}.
$$
\end{prop}

\subsection{Discontinuous time-process $Y$}
In several models the time-process $Y$ has discontinuous sample paths. Therefore, in this section
we investigate additional cases relaxing the assumption on continuity of sample paths of $Y$.

\begin{prop}\label{Thmfbm3}
Assume that $\theta>\alpha/2$, $c>0$  and the random variable $Y(T)$ possesses an absolutely continuous
distribution with probability density function 
which is regularly varying at $\infty$ with
index $\lambda+1>1$.
Then
\BQN
\lim_{u\to\infty}\frac{\log \psi_T(u)}{u^{\frac{2\theta-\alpha}{\theta}}}
= \frac{-1}{2V_0^{2}}.
\EQN
\end{prop}

\begin{prop}\label{Thmfbm4}
Assume that $\theta>\alpha/2$, $c>0$  and $Y(T)$ possesses absolutely continuous
distribution with probability density function $\rho_{Y(T)}(\cdot)$ such that
\tb{$\lim_{u\to\infty}\log (\rho_{Y(T)}(u))/u^p=-L$ for some $p,L>0$.}\\
i) If $2\theta-\alpha>p$, then
\BQN
\lim_{u\to\infty}\frac{\log \psi_T(u)}{u^{\frac{2\theta-\alpha}{\theta}}}
= \frac{-1}{2V_0^{2}}.
\EQN

ii) If $2\theta-\alpha<p$, then
\BQN
\lim_{u\to\infty}\frac{\log \psi_T(u)}
                      { u^{\frac{2p}{\alpha+p}} }
=
-\LT(\frac{1}{2}\LT(\frac{2p}{\alpha}\RT)^{\frac{\alpha}{\alpha+p}}+\LT(\frac{\alpha}{2p}\RT)^{\frac{p}{\alpha+p}}\RT)L^{\frac{\alpha}{\alpha+p}}.
\EQN

iii) If $2\theta-\alpha=p$, then
\BQN
\lim_{u\to\infty}\frac{\log \psi_T(u)}
{u^{\frac{2\theta-\alpha}{\theta}}}
=
-\LT(\frac{(1+c A_0^\theta)^2}{2A_0^\alpha}+LA_0^{2\theta-\alpha}\RT),
\EQN
where
$$
A_0=\left(\frac{c(\alpha-\theta)+\sqrt{c^2(\theta-\alpha)^2+2\alpha\LT(c^2(\theta-\alpha/2)+L(2\theta-\alpha)\RT)}}
{c^2(2\theta-\alpha)+2L(2\theta-\alpha)}\right)^{1/\theta}.
$$
\end{prop}
Proofs  of Propositions \ref{Thmfbm3} and \ref{Thmfbm4} are given in Section \ref{s.prop3_4}.

\begin{example}
Assume that $Y(t)=\sum_{i=1}^{N(t)} Z_i$ is a compound Poisson process with
$Z_i,i\ge 1$ non-negative independent random variables with  common 
probability density function $h(\cdot)$ which is regularly varying at $\infty$ with
index $\lambda+1>1$. If  $h(\cdot)$ is monotone and
$\{N(t),t \ge 0\}$ is a homogeneous Poisson process  with intensity $\mu>0$, then by Proposition \ref{Thmfbm3}
\[
\lim_{u\to\infty}\frac{\log \psi_T(u)}{u^{\frac{2\theta-\alpha}{\theta}}}
= \frac{-1}{2V_0^{2}}.
\]
\end{example}

\begin{example}
Let $\{\Gamma_t, t \geq 0\}$ be a Gamma process with parameter $\nu > 0$,
i.e., a L\'evy process such that the increments
\mbox{$\Gamma_{t+s} - \Gamma_t$} have Gamma  distribution with probability density function
\[
    f(x) = \frac{1}{\Gamma(\frac{s}{\nu})}
    x^{\frac{s}{\nu}-1}\exp(-x), \quad x>0.
\]
By fractional Laplace motion $\{ L_\alpha(t), t \geq 0 \}$ we denote a random process defined as
\[
    \{ L_\alpha(t), t \geq 0 \} \stackrel{d}{=} \{\sigma B_\alpha(\Gamma_t), t \geq 0 \}, \quad \sigma >0.
\]
A standard fractional Laplace motion corresponds to $\sigma = \nu
= 1$; see, e.g., Kozubowski et al. (2004), (2006).
Choosing $Y(t)=\Gamma_t$, we consider below finite-time ruin probability of risk process modelled by fractional Laplace motion
\BQNY
\psi_T(u)=\pk{\underset{t\in[0,T]}\inf \Bigl( u+c\Gamma_t-B_{\alpha}(\Gamma_t)\Bigr)<0}
=\pk{\underset{t\in[0,T]}\sup \Bigl(B_{\alpha}(\Gamma_t)-c\Gamma_t \Bigr)>u}.
\EQNY
For this model Proposition \ref{Thmfbm4} implies:

i) If $\alpha<1$, then
\BQN\nonumber
\lim_{u\to\infty}\frac{\log \psi_T(u)}{u^{2-\alpha}}
=
\frac{-2}{(2-\alpha)^2}\LT(\frac{c(2-\alpha)}{\alpha}\RT)^{\alpha}.
\EQN
ii) If $1<\alpha<2$, then
\BQN\nonumber
\lim_{u\to\infty}\frac{\log \psi_T(u)}
                      { u^{\frac{2}{\alpha+1}} }
=
-\frac{1}{2}\LT(\frac{2}{\alpha}\RT)^{\frac{\alpha}{\alpha+1}}-\LT(\frac{\alpha}{2}\RT)^{\frac{1}{\alpha+1}}.
\EQN

iii) If $\alpha=1$, then
\BQN\nonumber
\lim_{u\to\infty}\frac{\log \psi_T(u)}
{u}
= -\frac{(1+c A_0)^2}{2A_0}-A_0,
\EQN
where
$
A_0=\frac{1}
{\sqrt{c^2+2}}.
$
\end{example}

\COM{
\begin{remark}
\footnote{K: Maybe it suffices to have only the above examples?}
In this remark we present some candidates of  the time-changed process $\{Y(t), t\ge 0\}$.

1. $Y(t)=\Lambda t$, with $\Lambda$ a positive random variable satisfying \eqref{regV}, \eqref{eqwei} or \eqref{limL}.

2. $Y(t)=\sup_{s\in[0,t]}B_1(s)$ with $B_1$ a standard Brownian motion. In view of Lemma of .... xxx
\BQNY
\pk{Y(T)>u}=\pk{\sqrt{T}\abs{\mathcal{N}}>u}\ \sim \ \sqrt{\frac{2 T}{\pi}}u^{-1}\exp\LT(-\frac{u^2}{2T}\RT),\ \ u\to\IF,
\EQNY
implying that $Y(T)\in \mathcal{W}\LT(2, \frac{1}{2T}, \sqrt{\frac{2 T}{\pi}}, \xE{-1}\RT)$ and further $Y(T)$ has asymptotically a log-power tail (recall \eqref{limL}).

3. $Y(t)=\int_0^tB_1^2(s)ds$. it is known that (cf. Ex 4.1 of Vo{\ss} (2004))
\BQNY
\pk{Y(T)>u}\logsim \exp\LT(-\frac{\pi^2}{8T^2}u\RT)
\EQNY
showing that $Y(T)$ has asymptotically a log-power tail.
\end{remark}
}

\section{Proofs}
This section is dedicated to proofs of our results. In what follows, the positive constant $\const$ may be different from line to line.
We begin with a lemma which is of some interests on its own.

\begin{lem}\label{insensitive}
Let $X$ be a non-negative random variable which is $u^{1-p}$-insensitive, with $p>0$.
Then, for any positive constant  $B$
\[
\lim_{u\to\infty}\frac{\exp(-Bu^p)}{\pk{X>u}}=0.
\]
\end{lem}
\prooflem{insensitive} First observe that, 
\cJ{for any sufficiently large $u$ there exists some $\theta_u\in[0,1]$ such that}
(set next  $Y=X^{p}$)
\begin{eqnarray*}
\tn{
\pk{Y>u}}&\le&\pk{Y>u-1}\\
&=&
\pk{X>(u-1)^{1/p}}\\
& =&
\pk{X>u^{1/p}-(1/p)(u-\theta_u)^{1/p-1}}\\
&\le&
\pk{X>u^{1/p}-(2/p)(u)^{1/p-1}}.
\end{eqnarray*}
By insensitivity of $X$,
we immediately get that
\[
\pk{X>u^{1/p}-(2/p)(u)^{1/p-1}}
\sim
\pk{X>u^{1/p}}= \pk{Y>u}
\]
as $u\to\infty$.
Hence
\begin{eqnarray*}
 \pk{Y>u-1}\sim \pk{Y>u}, \quad u\to \IF. \label{1.ins}
\end{eqnarray*}
Consequently, for any  \tn{$\varepsilon\in(0,1)$} 
there exists $A>0$ such that for any
$v>A$
\tn{
we have
$\pk{Y>v}\ge(1-\varepsilon)\pk{Y>v-1}$. Thus, for sufficiently large $u$}
\begin{eqnarray}
\pk{Y>u}
&\ge&
(1-\varepsilon)\pk{Y>u-1}
\ge
(1-\varepsilon)^2\pk{Y>u-2}
\ge
(1-\varepsilon)^{u-A}\pk{Y>A}\nonumber
\end{eqnarray}
implying that for each $B>0$ 
\[
\lim_{u\rw\IF}\frac{\exp(-Bu)}{\pk{Y>u}}=0
\]
and thus the claim follows. \QED

\subsection{Proof of Theorem \ref{ThmT}}\label{s.proof.1}
i)
We first give the upper bound. For any $\ve\in (0,\tO), u>0$,  we derive that
\BQNY
\pk{\underset{t\in[0,\TT]}\sup X(t)>u}&\le&\pk{\underset{t\in[0, \TT ]}\sup X(t)>u, \ \TT< \tO-\ve}+\pk{\underset{ t\in[0,\TT] }\sup X(t)>u,\ \TT\ge \tO - \ve } \\
&\le&\pk{\underset{t\in[0, \tO - \ve ]}\sup X(t)>u}+\pk{\underset{ t\in[0,\infty) }\sup X(t)>u}\pk{\TT\ge \tO - \ve },
\EQNY
Similarly, for any $u>0$
\BQNY
\pk{\underset{t\in[0,\TT]}\sup X(t)>u}&\ge&\pk{\underset{t\in[0, \tO +\ve]}\sup X(t)>u}\pk{\TT\ge \tO + \ve }\\
 &\ge&\left(\pk{\underset{t\in[0, \IF)}\sup X(t)>u}-\pk{\underset{t\in[\tO +\ve,\IF)}\sup X(t)>u}\right)\pk{\TT\ge  \tO + \ve }.
\EQNY
Choosing $\ve$ small enough such that $\Lsup{t\to\IF}\ \sigma(t)<\sigma(\tO\pm\ve)<1$ and using
the Borell-TIS inequality (e.g., Adler and Taylor (2007)), we conclude that
\BQNY
\pk{\underset{t\in[0, \tO-\ve]}\sup X(t)>u}=\OU{\pk{\underset{t\in[0, \IF)}\sup X(t)>u}}\\
 \EQNY
and
 \BQNY
 \pk{\underset{t\in[\tO +\ve,\IF)}\sup X(t)>u}=\OU{\pk{\underset{t\in[0, \IF)}\sup X(t)>u}}
\EQNY
as $u\to\infty$.
Thus the claim of the first statement follows by letting $\ve\to 0$.\\
%
ii) Assume that $\TT$ is $u^{1-2/(\gamma(1+\beta))}$-insensitive and
let $\eU=u^{-2/(1+\beta)}$.
\cE{With} similar arguments as in the proof of statement $i)$  we obtain
\BQN
\pXT &=& \pk{\underset{t\in[0, \tn{\TT/u^\gamma }]}\sup X(t)>u}\notag \\
 &\le& \pk{\underset{t\in[0, \tO - \eU ]}\sup X(t)>u}+\pk{\underset{t\in[0,\infty) }\sup X(t)>u}\pk{\TT> ( \tO - \eU )u^{\gamma}}
\label{TUubound}
\EQN
and
\BQNY
\pXT&\ge& \pk{\underset{t\in[0, \tO + \eU ]}\sup X(t)>u}\pk{\TT> ( \tO + \eU )u^{\gamma}}\\
&\ge& \left(\pk{\underset{t\in [0,\infty) }\sup X(t)>u}-\pk{\underset{t\in[ \tO + \eU ,\infty)}\sup X(t)>u}\right)\pk{\TT> ( \tO + \eU )u^{\gamma}}.
\EQNY
Further, by insensitivity of $\TT$ as $u\to \IF$
\BQN
\pk{\TT> ( \tO\pm\eU )u^{\gamma}}&\sim& \pk{\TT> \tO u^{\gamma}}.\label{simT}
\EQN
\cE{Thanks to Piterbarg inequality} (cf.\ Theorem 8.1 of Piterbarg (1996)),
\tn{combined with {\bf A3},} \cE{for all $u$ large} we have
\BQNY
\pk{\underset{t\in[0, \tO - \eU ]}\sup X(t)>u}
&\le& \const u^{2/r-1}\exp\LT(-\frac{u^2}{2}\RT)\exp\LT(- \frac{a u^2( \eU )^\beta}{4}\RT)
\EQNY
for some positive constant $\const$ not depending on $u$.  Similarly,
using additionally Borell-TIS inequality, for $u$ large and some $G> t_0$
\BQNY
\pk{\underset{t\in[ \tO + \eU ,\infty)}\sup X(t)>u}&\le&\pk{\underset{t\in[ \tO + \eU ,G]}\sup X(t)>u}+\pk{\underset{t\in[ G,\infty)}\sup X(t)>u}\\
&\le& \const  u^{2/r-1} \exp\LT(-\frac{u^2}{2}\RT)\exp\LT(- \frac{a u^2( \eU )^\beta}{4}\RT)+\exp\LT(-\frac{\LT(u-\E{\underset{t\in[G,\IF)}\sup X(t)}\RT)^2}{2 \underset{t\in[G,\infty)}\sup \sigma^2(t)}\RT).
\EQNY
Hence, using that $\eU=u^{-2/(1+\beta)}$ and the results of \nelem{insensitive}
\BQNY
\lim_{u\rw\IF}\frac{\exp\LT(- \frac{a}{4} u^{\frac{2}{1+\beta}}\RT)}{\pk{\TT>t_0u^\gamma}}= 0 
\EQNY
and thus
\BQN
&&\pk{\underset{t\in[ \tO + \eU ,\infty)}\sup X(t)>u}=o\LT(\pk{\underset{ t\in[0,\infty) }\sup X(t)>u}\RT)\label{neg1},\\
&&\pk{\underset{t\in[0, \tO - \eU ]}\sup X(t)>u}=o\LT(\pk{\underset{ t\in[0,\infty) }\sup X(t)>u} \pk{\TT> ( \tO - \eU )u^{\gamma}} \RT).\label{neg2}
\EQN
Consequently, combining \eqref{TUubound}-\eqref{neg2} we obtain
\BQNY
\pXT&\sim& \pk{\underset{t\in [0,\infty) }\sup X(t)>u}\pk{\TT>  \tO u^\gamma }
\EQNY
as $u\to \IF$, which completes the proof. \QED

\subsection{ Proof of Theorem \ref{ThmLog}}\label{s.proof.2}
For the proof we need to distinguish between three different cases depending on the value of $\gamma p$.
\tn{For notational simplicity, let $\TU:=\TT/u^\gamma$.}\\
\underline{Case  $\gamma p< 2$}. For some  $\van >0$ small enough and $u$ large
\BQNY
\pXT&=&\pk{\underset{t\in[0,\TU]}\sup X(t)>u, \TU>u^{\frac{2-\gamma p}{p}-\van}}+\pk{\underset{t\in[0,\TU]}\sup X(t)>u, \TU\le u^{\frac{2-\gamma p}{p}-\van}}\\
&\ge&\pk{\TT>u^{\gamma+\frac{2-\gamma p}{p}-\van}}\pk{\underset{t\in\LT[0,u^{\frac{2-\gamma p}{p}-\van}\RT]}\sup X(t)>u}\\
&\ge&\pk{\TT>u^{\gamma+\frac{2-\gamma p}{p}-\van}}  \tn{\pk{X(t_0)>u}}\\
&\logsim& \pk{X(t_0)>u},
\EQNY
\tn{
which follows from the fact that, by the assumption on $\TT$,
$\pk{\TT>u^{\gamma+\frac{2-\gamma p}{p}-\van}}\logsim  \exp(-Lu^{2-\varepsilon p})$, while
$\pk{X(t_0)>u}\logsim \exp(-u^2/2)$.}
Since, in view of Theorem \ref{K1}
\BQNY
\pXT&\le&\pk{\underset{t\in[0,\infty)}\sup X(t)>u}\logsim \pk{X(t_0)>u},
\EQNY
then the claim \eqref{eqLogl} follows.

\underline{Case $\gamma p= 2$}.
\tn{Let $\widehat \sigma(s),\widetilde \sigma_{L,\gamma}(s), s\ge0$ and
$A_0$ be defined as in (\ref{not.1}) and (\ref{not.2}) respectively.}
The lower bound follows from the fact that
\BQN
\pXT&=&\pk{\underset{t\in[0,\TU]}\sup X(t)>u, \TU>A_0}+\pk{\underset{t \in[0,\TU]}\sup X(t)>u, \TU\le A_0}\nonumber\\
&\ge&\pk{\TT>A_0 u^\gamma}
\tn{
\pk{
\sigma(A_0)\mathcal{N} >u}}\nonumber\\
&\logsim&
\tn{
\exp\LT(-L A_0^{2/\gamma} u^2\RT)\exp\LT(-\frac{u^2}{2\sigma^2(A_0)}\RT)
}\nonumber\\
&=&\exp\LT(-{\cE{ \widetilde \sigma_{L,\gamma}(A_0)}}u^2\RT),
\label{clbound}
\EQN
where $\mathcal{N}$ is a standard Gaussian (i.e., \pE{an} $N(0,1)$) random variable.

Next we derive an upper bound. For some $0<m<1<M$ (to be determined later) we have
\BQNY
&&
\pXT=-\int_0^\IF\pk{\underset{t\in[0,s]}\sup X(t)>u}d\pk{\TU> s}\\
&&=-\int_0^{m A_0}\pk{\underset{t\in[0,s]}\sup X(t)>u}d\pk{\TU> s}-\int_{m A_0}^{M A_0}\pk{\underset{t\in[0,s]}\sup X(t)>u}d\pk{\TU> s}\\
&&-\int_{M A_0}^\IF\pk{\underset{t\in[0,s]}\sup X(t)>u}d\pk{\TU>s}\\
&&:=I_1(u)+I_2(u)+I_3(u).
\EQNY
We analyze $I_i(u), i=1,2,3$ separately. \\
Ad. $I_2(u)$.
Using Piterbarg inequality and integration by parts, we have that
\BQNY
I_2(u)&\le& -\int_{m A_0}^{M A_0}\LT(\const MA_0u^{2/r}\Psi\left(\frac{u}{\widehat \sigma(s) }\right)\RT)d\pk{\TU> s}\\
&=&-\const MA_0u^{2/r}\Bigg(\Psi\left(\frac{u}{\widehat \sigma(MA_0)}\right)\pk{\TT>MA_0 u^{\gamma}}-\Psi\left(\frac{u}{\widehat \sigma(mA_0)}\right)\pk{\TT>m A_0 u^{\gamma}}\\
&+&\int_{m A_0}^{MA_0} \frac{u}{\sqrt{2 \pi}}\LT(\frac{1}{\widehat \sigma(s) }\RT)'
\exp\Biggl(-\frac{u^2}{2\cE{\widehat \sigma^2(s)}}\Biggr)\pk{\TU> s} ds\Bigg)\\
&:=&I_{2,1}(u)+I_{2,2}(u)+I_{2,3}(u)\le I_{2,2}(u)+I_{2,3}(u).\label{eqI21}
\EQNY
Next we find bounds for $I_{2,i} (u), i=2,3$ one by one.
It straightforwardly follows that
\BQN\label{eqI22}
I_{2,2}(u)\logsim  \exp\Biggl(-\LT(\frac{1}{2 \widehat \sigma^2 (mA_0)}+L (mA_0)^{2/\gamma}\RT)u^2\Biggr)
\EQN
and, for any $\ve_2>0$ and $u$ large
\BQNY
I_{2,3}(u)&\le& \const MA_0\frac{u^{2/r+1}}{\sqrt{2 \pi}}\underset{s\in[mA_0,MA_0]}\max\LT(\frac{-1}{\widehat \sigma(s) }\RT)'
  \int_{m A_0}^{M A_0}\exp\LT(-\LT(\frac{1}{2\widehat \sigma^2(s)}+L(1-\ve_2)s^{2/\gamma}\RT)u^2\RT) ds\\
&\le& \const(MA_0)^{2}\frac{u^{2/r+1}}{\sqrt{2\pi}}\underset{s\in[mA_0,M A_0]}
 \max\LT(\frac{-1}{\widehat \sigma(s) }\RT)'
\exp\Biggl(-\underset{s\in[mA_0,MA_0]}\inf\LT(\frac{1}{2\widehat \sigma^2(s)}+L(1-\ve_2)s^{2/\gamma}\RT)u^2\Biggr).
\EQNY
Consequently,  letting $\ve_2\to 0$ we obtain
\BQN\label{eqI23}
\limsup_{u\to\infty}
\frac{\log I_{2,3}(u)}{u^2} \le - \widetilde \sigma_{L,\gamma}(A_0).
\EQN
Ad. $I_1(u)$.
By \pE{the} Piterbarg inequality we have that
\BQN\label{eqI1}
I_1(u)\le\pk{\underset{t\in[0,mA_0]}\sup X(t)>u}\le \const u^{2/r+1}\exp\LT(-\frac{u^2}{2\widehat \sigma^2(mA_0)}\RT),
\EQN
\tn{where $r\in(0,2]$ is as in {\bf A3}.}\\
Ad. $I_3(u)$.
We straightforwardly have that
\BQN\label{eqI3}
I_3(u)\le\pk{\TU>MA_0}\ \logsim\ \exp(-L(MA_0)^{2/\gamma}u^2).
\EQN
Now we are ready to determine both constants $m$ and $M$. First, choose $m$ such that
\BQNY
\frac{1}{2\widehat \sigma^2(mA_0)}>  \widetilde \sigma_{L,\gamma}(A_0),
\EQNY
 and then choose $M$ such that
\BQNY
L(MA_0)^{2/\gamma}> \widetilde \sigma_{L,\gamma}(A_0). 
\EQNY
We conclude from \eqref{eqI21}-\eqref{eqI3} that
\BQN\label{cubound}
\limsup_{u\to\infty}
\frac{\log
\pXT }
{u^2}
\le- \widetilde \sigma_{L,\gamma}(A_0). 
\EQN
Consequently, combination of (\ref{clbound}) with (\ref{cubound}) leads to
\BQNY
\lim_{u\to\infty}
\frac{\log
\pXT} 
{u^2}
=- \widetilde \sigma_{L,\gamma}(A_0).
\EQNY

\underline{Case $\gamma p> 2$ and $\sigma(0)>0$}.
\def\ue{u_\ve}
Let $\van\in(0,\frac{\gamma p-2}{p})$. Then we have (set $\ue := u^{\frac{2-\gamma p}{p}+\van}$)

\BQNY
\pXT&=&\pk{\underset{t\in[0,\TT_u]}\sup X(t)>u, \TT_u>\ue}+\pk{\underset{t\in[0,\TT_u]}\sup X(t)>u, \TU\le \ue}\\
&\le&\pk{\TT>u^{\frac{2}{p}+\van}}+\pk{\underset{t\in[0,\ue]}\sup X(t)>u}.
\EQNY
Since  $- \log \pk{\TT>u}$ is regularly varying at infinity with index $p$, then
\[
\lim_{u\to\infty}\frac{\log \pk{\TT>u^{\frac{2}{p}+\van}}}{u^2}=-\IF
\]
and, by Borell-TIS inequality,
for sufficiently large $u$
\[
\pk{\underset{t\in[0,\ue]}\sup X(t)>u}\le
2\exp\left(-\frac{\left(u-\mathbb{E}\left(\underset{t\in[0,\ue]} \sup X(t)\right)\right)^2}
                 {2\underset{t\in[0,\ue]} \max\sigma^2(t)}  \right).
\]
Combining the above with
\BQNY
\pXT&\ge& \pk{X(0)>u}
\EQNY
we obtain that
\[
\lim_{u\to\infty}\frac{\log\pXT}{u^2}=-\frac{1}{2\sigma^2(0)},
\]
which proves (\ref{eqLoggg}).\\

\underline{Case $\gamma p> 2$ and $\sigma(t)=Dt^\eta(1+\OU{1})$ as $t\downarrow0$}.
Let $g(t)=\frac{1}{2D^2t^{2\eta}}+Lt^p$, $t\ge0$,
which has a  unique minimum point at $t^*=\LT(\frac{\eta}{pLD^2}\RT)^{1/(2\eta+p)}$ with
$$
A_1:=g(t^*)=\frac{1}{2}D^{-\frac{2p}{2\eta+p}}(Lp/\eta)^{\frac{2\eta}{2\eta+p}}+L^{\frac{2\eta}{2\eta+p}}(\eta/(pD^2))^{\frac{p}{2\eta+p}}.
$$

Setting  $\mu=\frac{p\gamma-2}{2\eta+p}>0$ we may write
\BQNY
\pXT&=& \pk{\underset{t\in[0,\TU]}\sup X(t)>u, \TU>t^*u^{-\mu}}+\pk{\underset{t \in[0,\TU]}\sup X(t)>u, \TU\le t^*u^{-\mu}}\\
&\ge& \pk{\underset{t\in[0,t^*u^{-\mu}]}\sup X(t)>u}\pk{ \TT>t^*u^{\gamma-\mu}}\\
&\ge&  \pk{X(t^*u^{-\mu})>u}\pk{ \TT>t^*u^{\gamma-\mu}}\\
&\logsim&\exp\LT(-\frac{u^{2+2\eta\mu}}{2D^2(t^*)^{2\eta}}-L(t^*)^pu^{p(\gamma-\mu)}\RT)\\
&=&\exp\LT(-A_1u^{\frac{2\eta\gamma p+2p}{2\eta+p}}\RT).
\EQNY
In order to derive an upper bound, we replace $A_0$ with $t^*u^{-\mu}$ in the upper estimate of the case $\gamma p=2$.
Following step-by-step the same argument as in the upper bound of the case $\gamma p=2$,
we conclude that
\BQNY
\limsup_{u\to\infty}
\frac{\log \pXT} 
{u^{\frac{2\eta\gamma p+2p}{2\eta+p}}}
\le- A_1
\EQNY
and thus the proof is complete. \QED

\subsection{Proof of Propositions \ref{Thmfbm1} and \ref{Thmfbm2}}\label{s.prop1_2}
By the self-similar property of the fBm we see that
\BQN\label{eqfb1}
\pk{\underset{t\in[0,T]}\sup \Bigl(B_{\alpha}(Y(t))-c(Y(t))^\theta \Bigr)>u}& = &
\pk{\underset{t\in[0,Y(T)/u^{1/\theta}]}\sup \frac{B_{\alpha}(t)}{1+ct^\theta} >u^{1-\frac{\alpha}{2\theta}}}.
\EQN
The function $V(t)=\frac{t^{\alpha/2}}{1+ct^\theta}$ attains its maximum at the unique point
$$
\tO=\LT(\frac{\alpha/2}{c(\theta-\alpha/2)}\RT)^{1/\theta}
$$
and
$$
V_0=V(\tO)=\frac{2\theta-\alpha}{2\theta}\LT(\frac{\alpha}{c(2\theta-\alpha)}\RT)^{\frac{\alpha}{2\theta}}.
$$
Re-writing \eqref{eqfb1}, we have
\BQN\label{eqfb2}
\pk{\underset{t\in[0,T]}\sup \Bigl( B_{\alpha}(Y(t))-c(Y(t))^\theta \Bigr)>u}
& = &
\pk{\underset{t\in[0,Y(T)/u^{1/\theta}]}\sup Z(t) >V_0^{-1}u^{1-\frac{\alpha}{2\theta}}}\nonumber\\
&=&\pk{\underset{t\in\LT[0,\frac{Y(T)(V_0)^{-\frac{2}{2\theta-\alpha}}}{(v(u))^{\frac{2}{2\theta-\alpha}}}\RT]}\sup Z(t) >v(u)},
\EQN
where
$$
Z(t)=\frac{B_{\alpha}(t)}{t^{\alpha/2}}\frac{V(t)}{V_0},\ \ t\ge0\ \ \mathrm{and}\ \ v(u)=V_0^{-1}u^{1-\frac{\alpha}{2\theta}}.
$$
It is straightforward to check that
\BQNY
\sqrt{\E{Z^2(t)}}=1- \frac{1}{8}c^{\frac{2}{\theta}}\alpha^{1-\frac{2}{\theta}}(2\theta-\alpha)^{1+\frac{2}{\theta}}
(t-\tO)^2+\OU{(t-\tO)^2}
\EQNY
\tn{as $t\to t_0$}
and
\BQNY
\Co\LT(\frac{Z(t)}{\sqrt{\E{Z^2(t)}}}, \frac{Z(s)}{\sqrt{\E{Z^2(s)}}}\RT)=1- \frac{1}{2\tO^\alpha}
\abs{t-s}^\alpha+\OU{\abs{t-s}^\alpha}
\EQNY
as \tn{$s,t\to \tO$}.
Moreover, for any given positive constant $H$, it is derived that, for $t,s\in [0,H]$ with $\abs{t-s}<1$
\tn{
\BQNY
\E{Z(t)-Z(s)}^2
&=&
V_0^{-2}
\E{
\left(\frac{B_{\alpha}(t)}{1+ct^\theta}-\frac{B_{\alpha}(s)}{1+ct^\theta}\right)+
\left(\frac{B_{\alpha}(s)}{1+ct^\theta}-\frac{B_{\alpha}(s)}{1+cs^\theta}\right)}^2\nonumber\\
&\le&
2(V_0)^{-2}
\left(
\E{
\frac{B_{\alpha}(t)}{1+ct^\theta}-\frac{B_{\alpha}(s)}{1+ct^\theta}}^2
+
\E{
\frac{B_{\alpha}(s)}{1+ct^\theta}-\frac{B_{\alpha}(s)}{1+cs^\theta}}^2
\right)\nonumber\\
&=&
2(V_0)^{-2}
\left(
\frac{|t-s|^\alpha}{(1+ct^\theta)^2}
+
s^\alpha\LT(\frac{1}{1+ct^\theta}
-\frac{1}{1+cs^\theta}\RT)^2
\right)\\
&\le&\const \abs{t-s}^\alpha.
\EQNY
}
Thus assumptions {\bf A1-A4} are satisfied.
Additionally, by Theorem 1 in H\"{u}sler and  Piterbarg (1999)
\BQNY
\pk{\underset{t\in[0,\IF)}\sup Z(t) >v(u)}
\sim
Q\PH_\alpha
u^{\frac{(2\theta-\alpha)(2-\alpha)}{2\theta\alpha}}\Psi
\LT(\frac{2\theta}{2\theta-\alpha }\LT(\frac{\alpha}{c(2\theta-\alpha)}\RT)^{-\frac{\alpha}{2\theta}}u^{\frac{2\theta-\alpha}{2\theta}}\RT)
,\ \ \ u\to \IF,
\EQNY
where $Q:=2^{\frac{1}{2}+\frac{1}{\alpha}}\sqrt{\pi}c^{\frac{2-\alpha}{2\theta}}\alpha^{\frac{\alpha-2-\theta}{2\theta}}
\theta^{\frac{2-\alpha}{\alpha}}(2\theta-\alpha)^{\frac{\theta\alpha-4\theta+2\alpha-\alpha^2}{2\theta\alpha}}$.

The rest of the proof follows from an application of \netheo{ThmT} statement ii) and \netheo{ThmLog}. \QED

\subsection{Proof of Propositions \ref{Thmfbm3} and \ref{Thmfbm4}}\label{s.prop3_4}

Note that, for each $u>0$, using notation introduced in Section \ref{s.prop1_2}
\[
\pk{\underset{t\in[0,T]}\sup B_{\alpha}(Y(t))-c(Y(t))^\theta >u}
\le
\pk{\underset{t\in[0,Y(T)/u^{1/\theta}]}\sup Z(t) >V_0^{-1}u^{1-\frac{\alpha}{2\theta}}}.
\]
Thus it suffices to find logarithmically tight lower bounds for each subclass of densities of $Y(T)$.
The idea of the proof is the same both for the density of $Y(T)$ being regularly varying and
log-power tailed and heavily uses the idea of getting the lower bound in the proof of Theorem \ref{ThmLog}.
Hence we give only the argument for $Y(T)$ having regularly varying density function with index $\lambda+1$.
Under this scenario
\begin{eqnarray*}
\lefteqn{
\pk{\underset{t\in[0,T]}\sup \Bigl( B_{\alpha}(Y(t))-c(Y(t))^\theta \Bigr)>u}
\ge
\pk{\frac{B_{\alpha}(Y(T))}{u+c(Y(T))^\theta}>1}}\\
&\ge&
\min_{t\in[t_0u^{1/\theta}-u^{1/(2\theta)},t_0u^{1/\theta}+u^{1/(2\theta)}]}
\pk{\frac{B_{\alpha}(t)}{u+ct^\theta}>1}
\pk{Y(T)\in [t_0u^{1/\theta}-u^{1/(2\theta)},t_0u^{1/\theta}+u^{1/(2\theta)}]}.
\end{eqnarray*}
Using that
$
\pk{Y(T)\in [t_0u^{1/\theta}-u^{1/(2\theta)},t_0u^{1/\theta}+u^{1/(2\theta)}]}
$
is regularly varying at $\infty$
and
\[
\lim_{u\to\infty}
u^{\frac{\alpha}{\theta}-2} \log\LT(\min_{t\in[t_0u^{1/\theta}-u^{1/(2\theta)},t_0u^{1/\theta}+u^{1/(2\theta)}]}
\pk{\frac{B_{\alpha}(t)}{u+ct^\theta}>1}\RT) =
-\frac{1}{2V_0^2}
\]
we obtain a logarithmically tight lower bound, and thus the proof is complete. \QED
\\
\\
\tb{ {\bf Acknowledgement}: \pE{We are thankful to the referees for several comments and suggestions.}  K. D\c{e}bicki was partially supported
by NCN Grant No 2011/01/B/ST1/01521 (2011-2013). K. D\c{e}bicki,
E. Hashorva and L. Ji kindly acknowledge partial support by
the Swiss National Science Foundation Grant 200021-140633/1 and
by the project RARE -318984, a Marie Curie FP7 IRSES Fellowship. }

\end{document}